\documentclass[12pt]{article}

\usepackage{fullpage}
\usepackage{amsfonts}
\usepackage{latexsym,amssymb,graphics,epsfig}
\usepackage{amsmath, textcomp, subfigure, color}
\usepackage{multirow}
\usepackage{natbib}

\def\Figdir{./}

\newcommand{\ds}{\displaystyle}

\newcommand{\RR}{\mathbb{R}}
\newcommand{\NN}{\mathbb{N}}
\newcommand{\VV}{\mathbb{V}}

\def\UU{{\mathbb U}}
\def\XX{{\mathbb X}}
\def\SCS{\VV^0}

\def\obj{{\mathbb D}} 
\def\dynamics{f}
\def\one{y}
\def\two{z}
\def\controlONE{v}
\def\controlTWO{w}
\def\dynamicsTWO{R_{\two}}
\def\dynamicsONE{R_{\one}}
\def\catchONE{Y} 
\def\catchTWO{Z} 
\def\tonnes{\texttt{t}}

\newcommand{\ov}{\operatornamewithlimits}
\def\mtext#1{\,\mbox{#1}\,} 
\def\defegal{:=} 

\def\llower{^{\flat}}
\def\lloweropt{^{\flat,\star}}

\newtheorem{theorem}{Theorem}
\newtheorem{proposition}[theorem]{Proposition}

\newtheorem{corollary}[theorem]{Corollary}

\newenvironment{proof}{\small{\bf Proof.}}{\hfill$\Box$\normalsize
\bigskip}

\newcounter{Appctr}
\title{Viability Kernel for Ecosystem Management Models}

\author{
Eladio \textsc{Oca\~na Anaya} \footnote{\textsc{IMCA}, Instituto de Matem\'atica y Ciencias Afines,  Universidad Nacional de Ingenier\'{\i}a, Calle los Bi\'ologos 245, La Molina, Lima 12, Per\'u. eocana@imca.edu.pe, fax +511349-9838}
\and
Michel \textsc{De Lara}\footnote{%
Universit\'e Paris--Est, Cermics, 
 6--8 avenue Blaise Pascal, 77455 Marne la Vall\'ee Cedex 2, France. 
Corresponding author:  delara@cermics.enpc.fr, fax +33164153586}
\and 
Ricardo \textsc{Oliveros--Ramos} \thanks{\textsc{IMARPE}, 
Instituto del Mar del Per\'u,
Centro de Investigaciones en Modelado Oceanogr\'afico y Biol\'ogico Pesquero (CIMOBP), Apartado 22, Callao--Per\'u.
jtam@imarpe.gob.pe, roliveros@imarpe.gob.pe, fax +5114535053} 
\and
Jorge \textsc{Tam} \footnotemark[3]
}

\begin{document}

\maketitle

\begin{abstract}
We consider sustainable management issues formulated within the
framework of control theory.
The problem is one of controlling a discrete--time dynamical system
(\emph{e.~g.} population model) 
in the presence of state and control constraints, representing
conflicting economic and ecological issues for instance.
The viability kernel is known to play a basic role for the analysis of
such problems and the design of viable control feedbacks,
but its computation is not an easy task in general.
We study the viability of nonlinear generic ecosystem models 
under preservation and production constraints.
Under simple conditions on the growth rates at the boundary constraints,
we provide an explicit description of the viability kernel.
A numerical illustration is given for the hake--anchovy couple in the
Peruvian upwelling ecosystem. 
\medskip

\end{abstract}


\vspace{4 mm}

 \tableofcontents

\section{Introduction}
\label{sec:Introduction}

This paper deals with the  control of discrete--time dynamical systems 
of the form\\ $x(t+1)=\dynamics\big(x(t),u(t)\big)$, $t \in \NN$, 
with state $x(t) \in \XX$ and control $u(t) \in \UU$,
in the presence of state and control constraints $\big(x(t),u(t)\big)
\in \obj$.  The subset $\obj \subset \XX \times \UU$
describes ``acceptable configurations of the system''.
Such problems of dynamic control under constraints refer to viability 
\citep*{Aubin:1991} or invariance \citep*{Clarkeetal:1995} frameworks.
From the mathematical viewpoint, most of viability and weak invariance
results are addressed in the continuous time case. However, some
mathematical works deal with the discrete-time case. This  includes the
study of numerical schemes for the approximation of the viability
problems of the continuous dynamics as
in~\citep*{Aubin:1991,Stpierre:1994,QSP:1995}.  
In the control theory literature, problems of constrained control have
also been addressed in the discrete time case (see the survey paper
\citep*{Blanchini:1999});
reachability of target sets or tubes for nonlinear discrete time
dynamics is examined in~\citep*{bertsekas71minimax}. 

We consider sustainable management issues which can be formulated within
such a framework  as in
\citep*{benedoyen:2000,doyenbene2003,Bene-Doyen-Gabay:2001,Eisenack-Sheffran-Kropp:2006,martinetdoyenx,Mullonetal:2004,Rapaport-Terreaux-Doyen:2006,DeLara-Doyen-Guilbaud-Rochet_IJMS:2007,DeLara-Doyen:2008}.

The time index $t$ is an integer and the time period $[t,t+1[$ may be a
year, a month, etc. 
The dynamics is generally a population dynamics, with state vector $x(t)$
being either the biomass of a single species,
or a couple of biomasses for a predator--prey system,
or a vector of abundances at ages for one or for several species,
or abundances at different spatial patches, etc.
The control $u(t)$ may represent catches, harvesting mortality or harvesting effort.
The set $\obj$   may include biological, ecological and economic
objectives as in \citep*{Bene-Doyen-Gabay:2001}.
For instance, if the state $x$ is a vector of abundances at ages
and the control $u$ is a harvesting effort, 
$\obj = \{ (x,u) \mid B(x) \geq b\llower \; , E(x,u) \geq e\llower \}$ 
represents acceptable or acceptable configurations where 
conservation is ensured by an biological indicator $B(x) \geq b\llower$
(spawning stock biomass above a reference point, for instance)
and economics is taken into account via minimal catches 
$E(x,u) \geq e\llower$ (catches $E(x,u)$ above a threshold).

The viability kernel $\VV(\dynamics,\obj )$ associated with the dynamics
$\dynamics$ and the acceptable set $\obj$ is known to play a basic role
for the analysis of such problems and the design of viable control
feedbacks. Unfortunately, its computation is not an easy task in general.
In this paper, we focus on the case where the viability kernel
is given by a finite algorithm, with emphasis on population dynamics 
ecosystem models.
For nonlinear predator--prey systems, one can find 
descriptions of viability kernels 
for Lotka--Volterra systems 
in~\citep*{Bonneuil-Mullers:1997}
and population viability analysis in three
trophic-level food chains 
in~\citep*{Bonneuil-Saint-Pierre:2005}. 
In these latter references, the  control strategies are carried by the
species (``fitness-maximizing'' strategies), the constraint set
corresponds to having high enough densities, 
and the models are in continuous time.
In our approach, the controls are the harvesting efforts, the constraint
set includes both high abundances and high catches,
and the model is in discrete time.
\medskip

The paper is organized as follows.
Section~\ref{sec:Viability_and_sustainable_management} is devoted to
recalls on discrete--time viability and its 
possible use for sustainable management. 
For generic nonlinear ecosystem models, we
provide an explicit description of the viability kernel for preservation
and production constraints in
Section~\ref{sec:Viable_control_of_generic_nonlinear_ecosystem_models},
together with viable controls. 
An illustration in ecosystem management and numerical applications are
given for the hake--anchovy couple in the Peruvian upwelling ecosystem in 
Section~\ref{sec:A_numerical_application_to_the_hake-anchovy_Peruvian_ecosystem}.

\section{Discrete--time viability}
\label{sec:Viability_and_sustainable_management}

Let us  consider a nonlinear control system  described in discrete--time
by the difference equation
\begin{equation}
\left\{ \begin{array}{l}
x(t+1)=\dynamics\big(x(t),u(t)\big) \mtext{ for all } t \in \NN ,\\
x(0)=x_0 \mtext{ given,}
\end{array} \right .
\label{eq:generaldyn}
\end{equation}
where  the \emph{state variable} $x(t)$ belongs to  the finite
dimensional state space $\XX=\RR^{n_{\XX}}$, the \emph{control variable
} $u(t)$ is an element of the \emph{control set} $\UU=\RR^{n_{\UU}}$
while the \emph{dynamics} $\dynamics$ maps $\XX \times \UU$ into  $\XX$.

A controller or a decision maker  describes ``acceptable
configurations of the system'' through a set 
$\obj \subset \XX \times \UU$ termed the \emph{acceptable set}
\begin{equation}\label{eq:constraint}
\big(x(t),u(t)\big) \in \obj \mtext{ for all } t \in \NN \; ,
\end{equation}
where $\obj $ includes both system states and controls
constraints. 

The \emph{state constraints set} $\SCS$ associated with
$\obj $ is obtained by
projecting the acceptable set $\obj $ onto the state space $\XX$:
\begin{equation}
\SCS \defegal {\rm Proj}_{\XX}(\obj ) =
\{x \in \XX \mid \exists u \in \UU \, , \, (x,u) \in \obj \} \; .
\label{eq:state_constraints_set}
\end{equation}

Viability is defined as the ability to choose,
at each time step $t \in \NN$, a control $u(t) \in \UU$ such that the system
configuration remains acceptable.
More precisely, the system is viable  if the following feasible set is
not empty:
\begin{equation}
\VV(\dynamics,\obj ) \defegal \left\{x_0\in \XX \left|
\begin{array}{l}\exists\; (u(0),u(1), \ldots ) \mtext{ and } (x(0),x(1),
\ldots )\\
\mtext{ satisfying } \eqref{eq:generaldyn} \mtext{ and }\eqref
{eq:constraint}
\end{array} \right.
\right\} \; .
\label{eq:viability_kernel}
\end{equation}
The set $\VV(\dynamics,\obj )$ is called the \emph{viability kernel} 
\citep*{Aubin:1991} associated
with the dynamics $\dynamics$ and the acceptable set $\obj $.
By definition, we have $\VV(\dynamics,\obj) \subset \SCS
={\rm Proj}_{\XX}(\obj )$ but, in general, the
inclusion is strict. For a decision maker or control designer, knowing
the viability kernel has practical interest since it describes the
states from which controls can be found that maintain the
system in an acceptable configuration forever.
However, computing this kernel is not an easy task in general.
\bigskip

We now focus on the tools to achieve viability.
A subset $\VV$ is said to be
\emph{weakly invariant}  for the dynamics $\dynamics$ in the acceptable set
$\obj $, or a \emph{viability domain} of $\dynamics$ in $\obj $, if
\begin{equation}
\forall x \in \VV  \, , \quad \exists u \in \UU  \, , \quad
(x,u) \in \obj  \mtext{ and } \dynamics(x,u) \in \VV \; .
\label{eq:viability_domain}
\end{equation}
That is, if one starts from $\VV$, an acceptable control may
transfer the state in $\VV$.

Moreover, according to viability theory \citep*{Aubin:1991}, the viability
kernel $\VV(\dynamics,\obj )$ turns out to be the union of all viability
domains, or also the largest viability domain:
\begin{equation}
\VV(\dynamics,\obj)= \bigcup \biggl\{ \VV,\;\VV \subset \SCS,\;\VV\mtext{ 
viability domain for $\dynamics$ in ${\obj}$}  \biggr\}  \; .
\label{eq:viability_kernel_union}
\end{equation}
A major interest of such a property lies in the fact that
any viability domain for the dynamics $\dynamics$ in the acceptable set
$\obj $ provides a \emph{lower approximation} of the viability kernel.
An \emph{upper approximation} $\VV_k$ of the viability kernel
is given by the so called
\emph{viability kernel until time $k$
associated with $\dynamics$ in $\obj $}:
\begin{equation}
\VV_{k} \defegal \left\{x_0\in \XX \left|
\begin{array}{l}\exists\; (u(0),u(1), \ldots, u(k) ) \mtext{ and }
(x(0),x(1), \ldots, x(k) )\\
\mtext{ satisfying } \eqref{eq:generaldyn} \mtext{ for } t=0,\ldots,k-1 \\
\mtext{ and }\eqref{eq:constraint} \mtext{ for } t=0,\ldots,k
\end{array} \right.
\right\} \; .
\label{eq:Vtau}
\end{equation}
We have
\begin{equation}
\VV(\dynamics,\obj ) \subset \VV_{k+1}
\subset \VV_{k} \subset \VV_0 = \SCS \, \mtext{ for all } k \in \NN \; .
\label{eq:inclusionsVtau}
\end{equation}
It may be seen by induction that the decreasing sequence of
viability kernels until time $k$ satisfies
\begin{equation}
\VV_0 = \SCS \mtext{ and }
\VV_{k+1} = \left\{ x \in \VV_k \; \left| \; \exists u \in \UU \, ,
\, (x,u) \in \obj \mtext{ and } \dynamics(x,u) \in \VV_k
\right. \right\} \; .
\label{eq:induction}
\end{equation}
By~\eqref{eq:inclusionsVtau},
such an algorithm provides approximation from above of the
viability kernel as follows:
\begin{equation}
\VV(\dynamics,\obj ) \subset \ds \ov{\bigcap}_{k \in \NN} \VV_k =
\lim_{k \to +\infty} \! \downarrow  \! \VV_k \; .
\label{eq:algo1}
\end{equation}
Conditions ensuring that equality holds may be found in
\citep*{Stpierre:1994}. 
Notice that, when the decreasing sequence $(\VV_{k})_{k \in \NN}$ of 
viability kernels up to time $k$ is stationary, its limit is the 
viability kernel.
Indeed, if $\VV_{k} = \VV_{k+1} $ for some $k$, then $\VV_{k}$ is a
viability domain by~\eqref{eq:induction}.
Now, by~\eqref{eq:viability_kernel}, $\VV(\dynamics,\obj ) $ is the 
largest of viability domains. 
As a consequence, $\VV_{k} = \VV(\dynamics,\obj ) $ since
$ \VV(\dynamics,\obj ) \subset \VV_{k}$ by~\eqref{eq:inclusionsVtau}.
We shall use this property in the following
Sect.~\ref{sec:Viable_control_of_generic_nonlinear_ecosystem_models}. 
\bigskip

Once the viability kernel, or any approximation, or a viability
domain is known, we have to consider the management or control issue,
that is the  problem of selecting suitable controls at each time step.  
For any viability domain $\VV$ and any state $x \in \VV$, 
the following subset $\UU_{\VV}(x)$ of the decision set $\UU$ is not
empty:
\begin{equation}
\UU_{\VV}(x) \defegal
\{u \in \UU \mid (x,u) \in \obj  \mtext{ and } \dynamics(x,u) \in \VV \} \; .
\label{eq:viable-control}
\end{equation}
Therefore $\UU_{\VV(\dynamics,\obj)}(x)$ stands for the 
largest set of \emph{viable  controls associated with $x \in \XX$}.
Then, the decision design consists in the choice of
  a viable \emph{feedback} control,
namely any selection $\Psi: \; \XX \rightarrow \UU$
which associates with each state $x \in \VV(\dynamics,\obj)$ a control
$u=\Psi(x)$ satisfying $\Psi(x) \in \UU_{\VV(\dynamics,\obj)}(x)$.
\bigskip

In the context of sustainable management, viability concepts and
methods may help giving a framework for setting decision making.
First, one should delineate \emph{perpetual objectives} (the set
$\obj$). Second, \emph{operational objectives (advice)} can be
obtained as viable controls. This approach is illustrated in
\citep*{DeLara-Doyen-Guilbaud-Rochet_SCL:2006}, and especially in
\citep*{DeLara-Doyen-Guilbaud-Rochet_IJMS:2007} for fishery
management.

\section{Viable control of generic nonlinear ecosystem models}
\label{sec:Viable_control_of_generic_nonlinear_ecosystem_models}

For a generic 
ecosystem model,
we provide an explicit description of the viability kernel.
Then, we shall specify the results for predator--prey systems, in
particular for discrete-time Lotka--Volterra models.

\subsection{Dynamics and constraints for an ecosystem management model}

For simplicity, we consider a two--dimensional state model.
However, the following Proposition~\ref{pr:predator_prey}
may be easiliy extended to 
$n$--dimensional systems as long as  each species is harvested by a
specific device: one species, one harvesting effort.
The two--dimensional state vector $(\one,\two)$ represents biomasses  and the
two--dimensional control $(\controlONE,\controlTWO)$ is harvesting effort
of each species, 
respectively. The catches are thus $\controlONE \one$ and 
$\controlTWO \two$ (measured in
biomass).\footnote{%
In fact, any expression of the form $c(\one,\controlONE)$ would fit for
the catches in the following Proposition~\ref{pr:predator_prey} as
soon as $\controlONE \mapsto c(\one,\controlONE)$ 
is strictly increasing and goes
from $0$ to $+\infty$ when $\controlONE$ goes from $0$ to $+\infty$.
The same holds for $d(\two,\controlTWO)$ instead of $\controlTWO \two$.
}

The dynamics $\dynamics$ in~\eqref{eq:generaldyn} is given by
\begin{equation}
\dynamics(\one,\two,\controlONE,\controlTWO)=\left(\begin{array}{ccc}
\one\dynamicsONE(\one,\two,\controlONE)\\ 
\two\dynamicsTWO(\one,\two,\controlTWO)
\end{array}\right)~\mtext{ for all } (\one,\two,\controlONE,\controlTWO)\in \RR^4 \; ,
\label{eq:ecosystem}
\end{equation}
where $\dynamicsONE:\RR^3 \to \RR$ and $\dynamicsTWO:\RR^3 \to \RR$ are two
functions representing growth coefficients
(the growth rates being $\dynamicsONE -1$ and $\dynamicsTWO -1 $). 
We do not make assumptions on the signs of 
$\dynamicsONE$ and $\dynamicsTWO$; the restrictions on the domain of variation
of $(\one,\two,\controlONE,\controlTWO)$ will result from the requirement that the trajectories
belong to the acceptable set $\obj$, which will include signs considerations.
\medskip

The discrete--time control system~\eqref{eq:generaldyn} 
with dynamics $\dynamics$ given by~\eqref{eq:ecosystem} now writes:
\begin{equation}
\left\{ \begin{array}{rcl}
 \one(t+1) &=& \one(t) \dynamicsONE\big(\one(t),\two(t),\controlONE(t)\big) \\
 \two(t+1) &=& \two(t) \dynamicsTWO\big(\one(t),\two(t),\controlTWO(t)\big) \; .
\end{array} \right.
\end{equation}
The implicit assumption in such functional form~\eqref{eq:ecosystem} 
for the dynamics is
that, during one time period, the harvesting effort $\controlONE(t)$ of species $\one$
only affects the same species $\one(t+1)$ and not the other one $\two(t+1)$
(the reverse holds for $\controlTWO(t)$ and species $\two$). Of course, after two
periods, $\two(t+2)$ depends on $\one(t+1)$ 
which depends on $\controlONE(t)$ so that both efforts affect both species. 
Thus, the time period is assumed to be short enough for the impact on
one species of harvesting the other species to be negligible.
\medskip

The acceptable set $\obj$ is defined by 
\emph{minimal biomass thresholds} 
$\one\llower\geq 0$, $\two\llower\geq 0$
and \emph{minimal catch thresholds} 
$\catchONE\llower\geq 0$, $\catchTWO\llower\geq 0$:
\begin{equation}
\obj=\{\,(\one,\two,\controlONE,\controlTWO)\in \RR^4 \mid 
\one\geq \one\llower, \; 
\two\geq \two\llower,  \; 
\controlONE \one\geq\catchONE\llower,  \; 
\controlTWO \two\geq\catchTWO\llower\,\} \; .
\label{eq:thresholds}
\end{equation}
Such a set includes both \emph{preservation}
\begin{equation*}
 \one(t)\geq \one\llower \; , \quad \two(t)\geq \two\llower \; ,
\end{equation*}
and \emph{production} requirements
\begin{equation*}
 \controlONE(t)\one(t) \geq \catchONE\llower \; , \quad 
\controlTWO(t)\two(t) \geq \catchTWO\llower
\end{equation*}
when thresholds are positive.

\subsection{Expression of the viability kernel}

The following Proposition~\ref{pr:predator_prey} gives an explicit
description of the viability kernel, under some conditions on the
minimal thresholds. 

\begin{proposition}
Assume that the function $\dynamicsONE:\RR^3 \to \RR$ is 
continuously decreasing\footnote{%
In all that follows, a mapping $\varphi: \RR \to \RR$ is said to be
increasing if $x \geq x' \Rightarrow \varphi(x)  \geq \varphi(x')$. 
The reverse holds for decreasing. Thus, with
this definition, a constant mapping is both increasing and decreasing.
} in the control $\controlONE$  
and satisfies $\lim_{\controlONE \to +\infty}
\dynamicsONE(\one,\two,\controlONE)\leq 0$,
and that $\dynamicsTWO:\RR^3 \to \RR$ is 
continuously decreasing in the control variable $\controlTWO$,
and satisfies $\lim_{\controlTWO \to +\infty}
\dynamicsTWO(\one,\two,\controlTWO)\leq 0$. 
If the thresholds in~\eqref{eq:thresholds} are such that
the following growth coefficients are greater than one
\begin{equation}
\dynamicsONE(\one\llower,\two\llower,\frac{\catchONE\llower}{\one\llower})
\geq 1
\mtext{ and } 
\dynamicsTWO(\one\llower,\two\llower,\frac{\catchTWO\llower}{\two\llower})
\geq 1 \; ,
\label{eq:favorable_conditions}
\end{equation}
the viability kernel associated with the dynamics
$\dynamics$ in~\eqref{eq:ecosystem} 
and the acceptable set $\obj$ in~\eqref{eq:thresholds} is given by
\begin{equation}
\VV(\dynamics,\obj) = \left\{(\one,\two) \mid
\one\geq \one\llower, \;  \two\geq \two\llower, \;  
\one\dynamicsONE(\one,\two,\frac{\catchONE\llower}{\one})\geq \one\llower, \;
\two\dynamicsTWO(\one,\two,\frac{\catchTWO\llower}{\two})\geq \two\llower 
\right\}.
\label{eq:predator_prey_viability_kernel}
\end{equation}
\label{pr:predator_prey}
\end{proposition}

Before giving the proof, let us comment the assumptions.
That the growth coefficients are decreasing with respect to the
harvesting effort is a natural assumption.
Conditions~\eqref{eq:favorable_conditions} mean that, at the point 
$(\one\llower,\two\llower)$ and applying efforts 
$u\llower=\frac{\catchONE\llower}{\one\llower}$,
$v\llower=\frac{\catchTWO\llower}{\two\llower}$,
the growth coefficients are greater than one, hence both populations grow;
hence, it could be thought that computing the viability kernel is
useless since everything looks fine. 
However, if all is fine at the point $(\one\llower,\two\llower)$, 
it is not obvious that this also goes for a larger domain.
Indeed, the ecosystem dynamics $\dynamics$ given by~\eqref{eq:ecosystem}
has no  monotonocity properties that would allow to extend a result 
valid for a point to a whole domain.
What is more, if continuous-time viability results mostly relies upon
assumptions at the frontier of the constraints set, this is no longer
true for discrete-time viability.
\medskip


\begin{proof}
According to induction~\eqref{eq:induction}, we have:
\begin{eqnarray*}
\VV_0 &=& \{\,(\one,\two) \left|  \one\geq \one\llower,  \two\geq \two\llower
\right. \,\},\\[3mm]
\VV_1 &=& \left\{(\one,\two) \left| 
\begin{array}{ccc}
\one\geq \one\llower,  \two\geq \two\llower 
\mtext{ and, for some }\, (\controlONE,\controlTWO)\geq 0,\\ 
\controlONE \one\geq \catchONE\llower, \controlTWO \two\geq \catchTWO\llower, 
\one\dynamicsONE(\one,\two,\controlONE)\geq \one\llower,
\two\dynamicsTWO(\one,\two,\controlTWO)\geq \two\llower
\end{array} \right. \right\} \\
&=& \left\{(\one,\two) \left| 
\one\geq \one\llower,  \two\geq \two\llower,  
\one\dynamicsONE(\one,\two,\frac{\catchONE\llower}{\one})\geq \one\llower,
\two\dynamicsTWO(\one,\two,\frac{\catchTWO\llower}{\two})\geq \two\llower \right. \right\}\\
&& \mtext{ because } \controlONE \mapsto \dynamicsONE(\one,\two,\controlONE) \mtext{ and } 
 \controlTWO \mapsto \dynamicsTWO(\one,\two,\controlTWO) \mtext{ are decreasing,}\\
&& \mtext{ and thus we may select } 
\controlONE=\frac{\catchONE\llower}{\one}, \; 
 \controlTWO =\frac{\catchTWO\llower}{\two}. \\
&& \mtext{Denoting } \one'=
\one\dynamicsONE(\one,\two,\controlONE),
\,\, \two'=\two\dynamicsTWO(\one,\two,\controlTWO), \mtext{we obtain,}\\[4mm]
\VV_2 &=& \left\{(\one,\two) \left| 
\begin{array}{l}
\one\geq \one\llower,  \two\geq \two\llower\mtext{ and, for some }\, (\controlONE,\controlTWO)\geq 0,\\ 
\controlONE \one\geq \catchONE\llower,\,\, \controlTWO \two\geq \catchTWO\llower\\
\one'\geq \one\llower,\,\, 
\one'\dynamicsONE(\one',\two',\frac{\catchONE\llower}{\one'})\geq \one\llower,\,\, 
\two'\geq \two\llower,\,\, 
\two'\dynamicsTWO(\one',\two',\frac{\catchTWO\llower}{\two'})\geq \two\llower
\end{array} \right. \right\} \; .
\end{eqnarray*}
We shall now make use of the property, recalled in
Sect.~\ref{sec:Viability_and_sustainable_management}, that
when the decreasing sequence $(\VV_{k})_{k \in \NN}$ of 
viability kernels up to time $k$ is stationary, its limit is the 
viability kernel $ \VV(\dynamics,\obj )$.
Hence, it suffices to show that $\VV_1\subset \VV_2$ to obtain that
$ \VV(\dynamics,\obj) = \VV_1$.
Let $(\one,\two)\in \VV_1$, so that
\begin{equation*}
\one\geq \one\llower, \quad \two\geq \two\llower \mtext{ and }~ 
\one\dynamicsONE(\one,\two,\frac{\catchONE\llower}{\one})\geq \one\llower, \quad 
\two\dynamicsTWO(\one,\two,\frac{\catchTWO\llower}{\two})\geq \two\llower \; .
\end{equation*}
Since $\dynamicsONE:\RR^3 \to \RR$ is
continuously decreasing in the control variable, with
$\lim_{\controlONE \to +\infty}\dynamicsONE(\one,\two,\controlONE)\leq 0$, 
and since $\one\dynamicsONE(\one,\two,\frac{\catchONE\llower}{\one})
\geq \one\llower$, there exists a $\hat{\controlONE} \geq
\frac{\catchONE\llower}{\one}$ (depending on $\one$ and $\two$) such that
$\one'=\one\dynamicsONE(\one,\two,\hat{\controlONE})= \one\llower$.
The same holds for $\dynamicsTWO:\RR^3 \to \RR$ and
$\two'=\two\dynamicsTWO(\one,\two,\hat{\controlTWO})= \two\llower$.
By~\eqref{eq:favorable_conditions}, we deduce that
\[
\one'\dynamicsONE(\one',\two',\frac{\catchONE\llower}{\one'}) = \one\llower
\dynamicsONE(\one\llower,\two\llower,\frac{\catchONE\llower}{\one\llower})
\geq \one\llower \mtext{ and }
\two'\dynamicsTWO(\one',\two',\frac{\catchTWO\llower}{\two'}) = \two\llower
\dynamicsTWO(\one\llower,\two\llower,\frac{\catchTWO\llower}{\two\llower})
\geq \two\llower \; .
\]
The inclusion $\VV_1\subset \VV_2$ follows.
\end{proof}

As a direct consequence of the proof of
Proposition~\ref{pr:predator_prey}, for each
$(\one,\two)\in\VV(\dynamics,\obj)$, the control $(\hat{\controlONE},\hat{\controlTWO})$ 
belongs to $\UU_{\VV(\dynamics,\obj)}(\one,\two)$, the set of viable controls
defined in \eqref{eq:viable-control}. 
More explicitly, we have the following expression of the viable controls
set, which results from the above observation and proof.

\begin{corollary}
Suppose that the assumptions of Proposition~\ref{pr:predator_prey} are
satisfied. 
Denoting 
\[
\left\{ \begin{array}{rcl}
\hat{\controlONE}(\one,\two) 
&=& \max \{ \controlONE \geq \displaystyle\frac{\catchONE\llower}{\one} \mid
\one\dynamicsONE(\one,\two,\controlONE)= \one\llower \} \\[3mm]
\hat{\controlTWO}(\one,\two) &=& \max \{ \controlTWO \geq \displaystyle\frac{\catchTWO\llower}{\two} \mid
\two\dynamicsTWO(\one,\two,\controlTWO)=\two\llower \} \; ,
\end{array} \right.
\]
the set of viable controls is given by
\begin{equation*}
\UU_{\VV(\dynamics,\obj)}(\one,\two)= \left\{(\controlONE,\controlTWO) \left| 
\begin{array}{l}
\hat{\controlONE}(\one,\two) \geq \controlONE\geq \displaystyle\frac{\catchONE\llower}{\one}, \quad 
\hat{\controlTWO}(\one,\two) \geq \controlTWO\geq \displaystyle\frac{\catchTWO\llower}{\two}, \\

\one'\dynamicsONE(\one',\two',\displaystyle\frac{\catchONE\llower}{\one'})\geq \one\llower,\,\, 
\two'\dynamicsTWO(\one',\two',\displaystyle\frac{\catchTWO\llower}{\two'})\geq \two\llower
\end{array} \right. \right\} 
\end{equation*}
where 
$\one'= \one\dynamicsONE(\one,\two,\controlONE)$,
$\two'=\two\dynamicsTWO(\one,\two,\controlTWO) $.
\label{cor:viable-contr}
\end{corollary}

\subsection{Viability kernel for a nonlinear predator--prey system}

Till now, no trophic relationship has been specified between species
$\one$ and $\two$. 
From now on, we shall focus on nonlinear predator--prey systems.
The two--dimensional state vector $(\one,\two)$ represents
biomasses of preys ($\one$) and predators ($\two$). 
This interpretation results from the monotonicity
assumptions made on function $\dynamicsTWO$ in the following
Proposition~\ref{pr:predator_prey_withoud_dd}, which
provides a simpler version of Proposition~\ref{pr:predator_prey}
when the predator $\two$ does not exhibit density--dependence. 

\begin{proposition}
Suppose that the dynamics $\dynamics$ is given by
\begin{equation}
\dynamics(\one,\two,\controlONE,\controlTWO)=\left(\begin{array}{ccc}
\one\dynamicsONE(\one,\two,\controlONE)\\ 
\two\dynamicsTWO(\one,\controlTWO)
\end{array}\right)~\mtext{ for all } 
(\one,\two,\controlONE,\controlTWO)\in \RR^4 \; .
\label{eq:predator_prey}
\end{equation}
Assume that the function $\dynamicsONE:\RR^3 \to \RR$ is decreasing in the
control $\controlONE$ and $\dynamicsTWO:\RR^2 \to \RR$ is increasing in the 
state variable $\one$ and continuously decreasing in the control
variable $\controlTWO$, and satisfies 
$\lim_{\controlTWO \to +\infty}\dynamicsTWO(\one,\controlTWO)\leq 0$. If 
\begin{equation}
\dynamicsONE(\one\llower,\two\llower,\frac{\catchONE\llower}{\one\llower})
\geq 1 \mtext{ and } 
\dynamicsTWO(\one\llower,\frac{\catchTWO\llower}{\two\llower})\geq 1 \; ,
\label{eq:favorable_conditions_withoud_dd}
\end{equation}
the viability kernel associated with the dynamics
$\dynamics$ in~\eqref{eq:predator_prey}
and the acceptable set $\obj$ in~\eqref{eq:thresholds} is given by
\begin{equation}
\VV(\dynamics,\obj) = \left\{(\one,\two) \mid
\one\geq \one\llower, \;  \two\geq \two\llower, \;  
\one\dynamicsONE(\one,\two,\frac{\catchONE\llower}{\one})\geq \one\llower
\right\}.
\label{eq:predator_prey_viability_kernel_withoud_dd}
\end{equation}
\label{pr:predator_prey_withoud_dd}
\end{proposition}

\begin{proof}
The proof follows the same lines as the one of
Proposition~\ref{pr:predator_prey}.
The description of the viability kernel is simpler because of the
following argument.
The monotonicity assumptions on each component of
function $\dynamicsTWO$ and condition
$\;\dynamicsTWO(\one\llower,\frac{\catchTWO\llower}{\two\llower})\geq 1$ 
lead to 
\[
\one\geq \one\llower \mtext{ and } \two\geq \two\llower \Rightarrow 
\dynamicsTWO(\one,\frac{\catchTWO\llower}{\two})\geq
\dynamicsTWO(\one,\frac{\catchTWO\llower}{\two\llower})\geq
\dynamicsTWO(\one\llower,\frac{\catchTWO\llower}{\two\llower})\geq
1=\frac{\two\llower}{\two\llower}\geq \frac{\two\llower}{\two} \; .
\]
Therefore, $\two\dynamicsTWO(\one,\frac{\catchTWO\llower}{\two})\geq \two\llower$
for all $\one\geq \one\llower \mtext{ and } \two\geq \two\llower$.
This explains why \eqref{eq:predator_prey_viability_kernel}
reduces to 
\eqref{eq:predator_prey_viability_kernel_withoud_dd}.
\end{proof}

\subsection{Viability kernel for a Lotka--Volterra system}

Consider the following discrete--time Lotka--Volterra system of
equations with density--dependence in the prey
\begin{equation}
\left\{ \begin{array}{rcl}
\one(t+1) &=& \displaystyle 
R\one(t)-\frac{R}{\kappa}\one^2(t)-\alpha \one(t)\two(t) 
- \controlONE(t) \one(t) \; ,\\[3mm]
\two(t+1) &=& L\two(t)+\beta \one(t)\two(t)-\controlTWO(t)\two(t) \; ,
\end{array} \right.
\label{eq:LV_with_dd}
\end{equation}
where $R>1$, $0<L<1$, $\alpha >0$, $\beta >0$ and
$\kappa=\tfrac{R}{R-1}K$, with $K>0$ the carrying capacity for prey. 
The dynamics~$\dynamics$ is as in~\eqref{eq:predator_prey},
with $\dynamicsONE(\one,\two,\controlONE)=
R-\frac{R}{\kappa}\one-\alpha \two - \controlONE $ and
$\dynamicsTWO(\one,\controlTWO)=L+\beta \one- \controlTWO $.
Proposition~\ref{pr:predator_prey_withoud_dd} gives the following
Corollary.

\begin{corollary}
Consider the Lotka--Volterra predator--prey model~\eqref{eq:LV_with_dd}.
Whenever 
\begin{equation*}
\one\llower\geq 
\frac{1-L}{\beta}
\mtext{ and }  
\two\llower\leq 
\frac{R-1}{\alpha} - \frac{R(1-L)}{\alpha\beta\kappa} \; ,
\end{equation*}
any minimal catch thresholds $\catchONE\llower$ and $\catchTWO\llower$
such that 
\begin{subequations}
\begin{align}
\catchONE\llower \leq & \catchONE\lloweropt \defegal
\one\llower(R-\tfrac{R}{\kappa}\one\llower -\alpha \two\llower -1) 
\label{eq:maximal_catches_LV_with_ddONE}
 \\[3mm]
\catchTWO\llower \leq &  
\catchTWO\lloweropt \defegal \two\llower(L + \beta \one\llower - 1) 
\label{eq:maximal_catches_LV_with_ddTWO}
\; , 
\end{align}
\end{subequations}
satisfy~\eqref{eq:favorable_conditions_withoud_dd} and
the viability kernel associated with the dynamics
$\dynamics$ in~\eqref{eq:LV_with_dd}
and the acceptable set $\obj$ in~\eqref{eq:thresholds} is given by 
\begin{equation}
\VV(\dynamics,\obj) = 
\left\{(\one,\two) \left|
\begin{array}{c}
\displaystyle 
\one\geq \one\llower, \; \two\llower \leq
\two\leq \frac{1}{\alpha} \left\lbrack R  \left( \frac{\kappa-\one}{\kappa} \right)- 
\frac{\catchONE\llower + \one\llower}{\one}\right\rbrack
\end{array}
\right. \right\} \; .
\label{eq:LV_kernel_with_dd}
\end{equation}
\label{cor:LV_kernel_with_dd}
\end{corollary}

This Corollary has the following practical consequence.
An initial state $(\one,\two) $ such that 
$\one\geq \one\llower$ and $\two \geq \two\llower$ belongs to 
all the viability kernels $\VV(\dynamics,\obj)$ associated to
$\catchTWO\llower \leq \catchTWO\lloweropt $, 
by~\eqref{eq:maximal_catches_LV_with_ddTWO},
and to
$ \catchONE\llower \leq \min \left\{ \catchONE\lloweropt ,
\one ( R - R \one/\kappa -   \alpha\two) -\one\llower  \right\}$, 
by~\eqref{eq:maximal_catches_LV_with_ddONE} and
\eqref{eq:LV_kernel_with_dd},
as long as this last quantity is nonnegative.
In other words, if viably managed, the fishery could produce at least 
$ \min \left\{ \catchONE\lloweropt ,
\one ( R - R \one/\kappa -   \alpha\two) -\one\llower  \right\} $ 
and $\catchTWO\lloweropt $, while respecting
biological thresholds $\one\llower$ and $\two\llower$.
We shall use this property in the following numerical application.

\section{Numerical application to the hake--anchovy couple in the 
Peruvian upwelling ecosystem}
\label{sec:A_numerical_application_to_the_hake-anchovy_Peruvian_ecosystem}

We provide a viability analysis of the
hake--anchovy Peruvian fisheries between the years 1971 and 1981.
For this, we shall consider a discrete-time Lotka--Volterra model 
for the couple anchovy (prey $\one$) and hake (predator $\two$), for
which the viability kernel has explicit description. 
The emphasis is not on developing a biological model,
but rather on decision-making using such a model. 

\begin{figure}
\centering
\begin{tabular}{cc}
\subfigure[Anchovy]{ \includegraphics[width=0.45\textwidth]%
{\Figdir 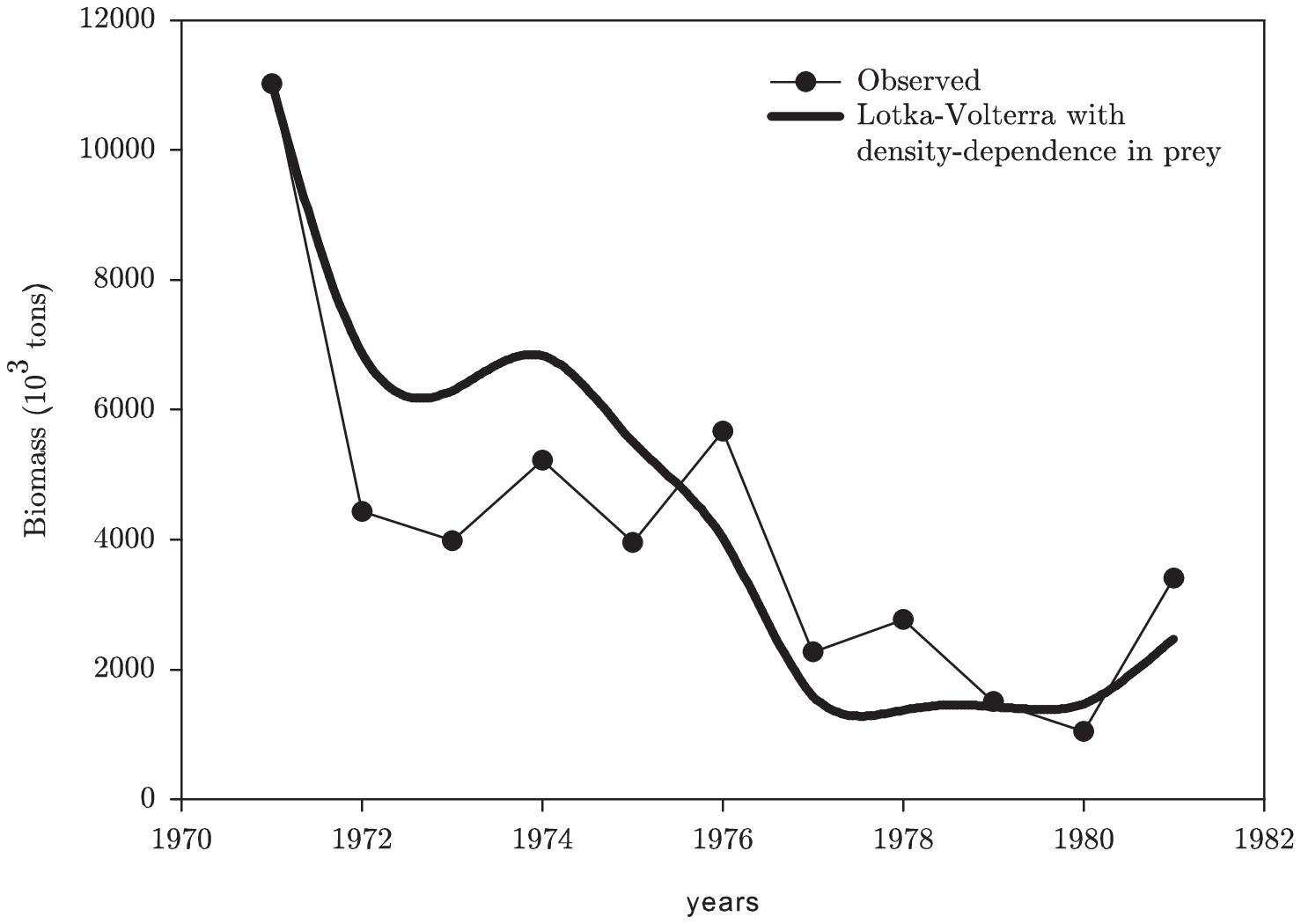}} &
\subfigure[Hake]{\includegraphics[width=0.45\textwidth]%
{\Figdir 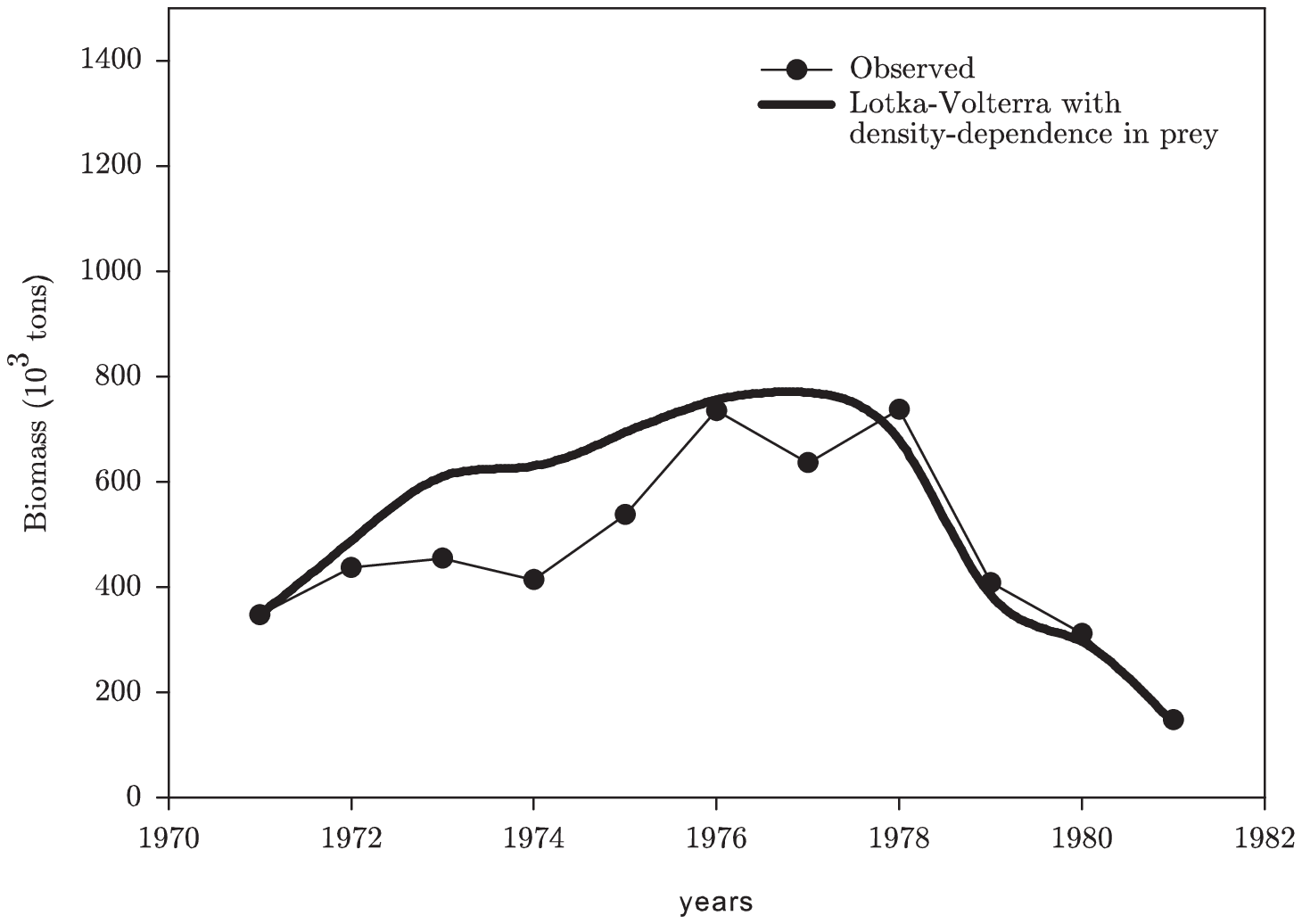}}
\end{tabular}
\caption{Comparison of observed and simulated biomasses of anchovy and
  hake using a Lotka--Volterra model with density-dependence in the
  prey. Model parameters are 
$R=2.25$~year$^{-1}$, 
$L=0.945$~year$^{-1}$, 
$\kappa=67~113~\times 10^3$~\tonnes\ ($K=37~285~\times 10^3$~\tonnes),
$\alpha=1.220\times 10^{-6}$~\tonnes$^{-1}$,
$\beta= 4.845\times 10^{-8}$~\tonnes$^{-1}$. }
\label{fig:anchovy_hake_non_linear.eps}
\end{figure}

The period between the years 1971 and 1981 has been chosen
because the competition between the fishery and hake was reduced due to
low anchovy catches, and because of the absence a strong warm event after El Ni\~{n}o 1972. We have 11 couples of biomasses, and the same for catches.
The 5 parameters of the model are estimated minimizing a weighted
residual squares sum function using a conjugate gradient method, with
central derivatives.
Estimated parameters and comparisons of observed and simulated biomasses are
shown in Figure~\ref{fig:anchovy_hake_non_linear.eps}.  

We consider values of $\one\llower=7~000~000$~\tonnes\ and
$\two\llower=200~000$~\tonnes\ for minimal biomass thresholds and values of
$\catchONE\llower=2~000~000$~\tonnes\ and 
  $\catchTWO\llower=5~000$~\tonnes\ for minimal catch thresholds 
\citep*{IMARPE:2000, IMARPE:2004}. 
Conditions~\eqref{eq:favorable_conditions_withoud_dd} in
Proposition~\ref{pr:predator_prey_withoud_dd} are satisfied 
with these values.
Indeed, the expressions 
in~\eqref{eq:maximal_catches_LV_with_ddONE}--\eqref{eq:maximal_catches_LV_with_ddTWO}
give:
\begin{equation}
\left\{
\begin{array}{rclclcl}
\catchONE\llower &= & 2~000~000~\tonnes & \leq & \catchONE\lloweropt &= &
5~399~000~\tonnes \\[3mm]
\catchTWO\llower  &= & 5~000~\tonnes & \leq &  
\catchTWO\lloweropt &=& 56~800~\tonnes \; . 
\end{array}
\right. 
\label{eq:numerical_maximal_catches_LV_with_dd} 
\end{equation}


\begin{figure}
	\centering
\includegraphics[width=0.8\textwidth]%
{\Figdir 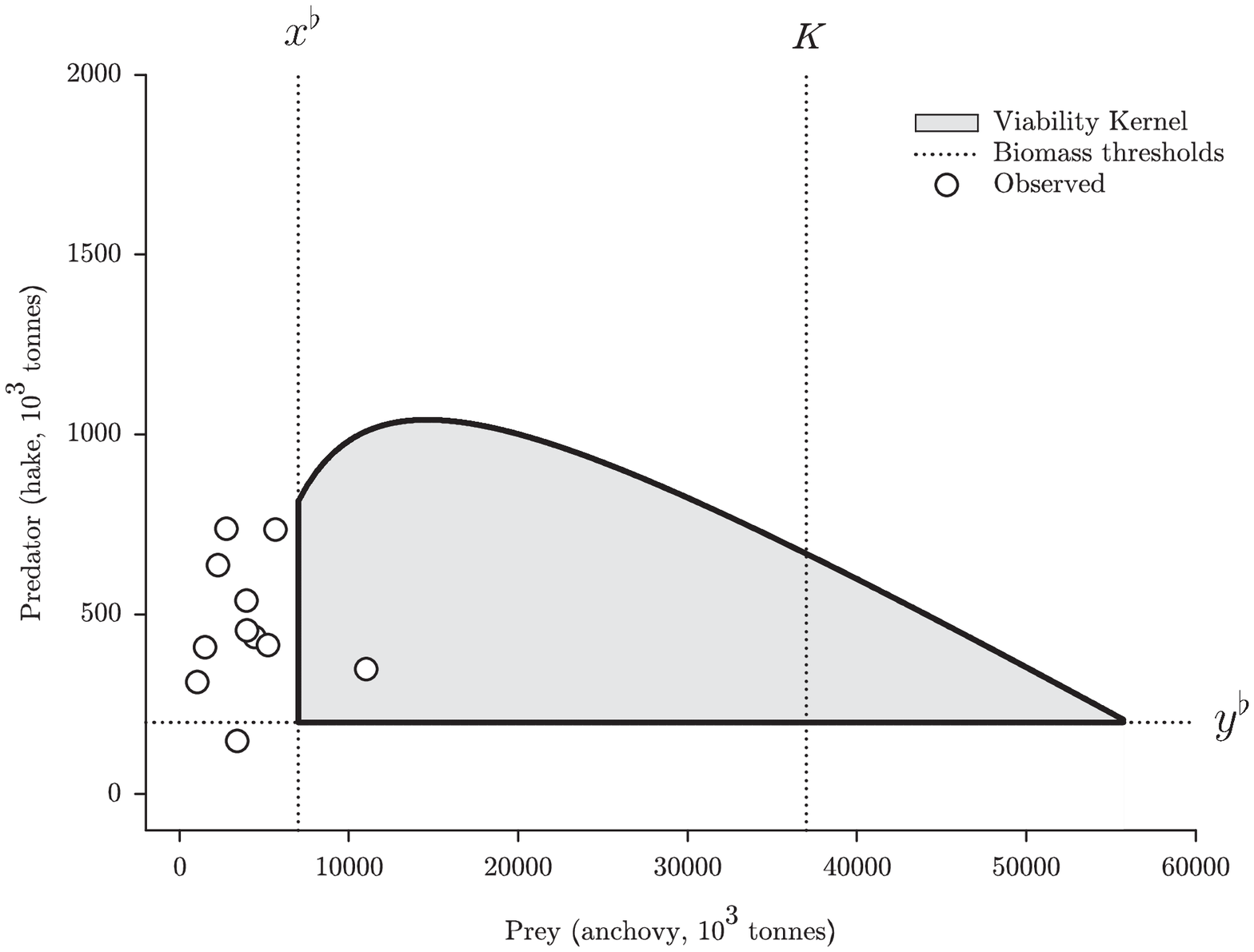}
\caption{Viability kernel (in grey) for a Lotka--Volterra model with
  density-dependence in the prey in the predator--prey phase space 
(with {$\one\llower=7\,000\,000$~\tonnes,
  $\two\llower=200\,000$~\tonnes,
  $\catchONE\llower=2\,000\,000$~\tonnes,
  $\catchTWO\llower=5\,000$~\tonnes}).
The unique point within the viability kernel is the initial point, while 
all subsequent points in the trajectory are outside the state constraint
set.
}
\label{fig:viab_kernel_non_linear.eps}
\end{figure}

The viability kernel is depicted in 
Figure~\ref{fig:viab_kernel_non_linear.eps}. 
The unique viable point (within the viability kernel) is the initial
point. 
Thus, based upon this model, the fishery could have been managed 
-- with appropriate viable controls -- 
to produce catches above $\catchONE\llower=2~000~000$~\tonnes\ and 
  $\catchTWO\llower=5~000$~\tonnes,
while ensuring biological conservation.
What is more, due to the remark following
Corollary~\ref{cor:LV_kernel_with_dd}, captures up to 
$\catchONE\llower=5~399~000$~\tonnes\ and 
  $\catchTWO\llower=56~800$~\tonnes\ were theoreticaly achievable in a
  sustainable way starting from year 1971.

\section{Conclusion}

Motivated by viable management of ecosystems, we have provided a general
condition ensuring an explicit construction of viability kernels and
have applied this to the viability analysis of generic ecosystem models
with harvesting. 

Our results have then been applied to a Lotka--Volterra model using the 
anchovy--hake couple in the Peruvian upwelling ecosystem.
We showed that, during the anchovy collapse, 
theoretically the fishery could have been viably managed  
to produce catches above the expected threshold levels 
while ensuring biological conservation.

It is interesting to notice that the kind of 
maximum sustainable yields 
$\catchONE\lloweropt = 5~399~000~\tonnes $ and
$ \catchTWO\lloweropt = 56~800~\tonnes $ 
provided by our approach 
in~\eqref{eq:numerical_maximal_catches_LV_with_dd} 
are comparable to the $4~250~000$~\tonnes\ anchovy yield and the 
$55~000$~\tonnes\ hake yield, respectively, established for the year 2006
\citep*{PRODUCE:2005, PRODUCE:2006b, PRODUCE:2006c},
or to the $5~000~000$~\tonnes\ anchovy yield 
and the $35~000$~\tonnes\ hake yield,
respectively, established for the year 2007 \citep*{PRODUCE:2006a, PRODUCE:2007a, PRODUCE:2007b}. 
So, despite simplicity\footnote{%
In addition to hake, there are other important
predators of anchovy in the Peruvian upwelling ecosystem, such as
mackerel and horse mackerel, seabirds and pinnipeds, which were not
considered. Also,
anchovy has been an important prey of hake, but other prey species have
been found in the opportunistic diet of hake \citep*{Tam:2006}
}
of the models considered, our approach may provide a mean of designing
sustainable yields from an ecosystem point of view. 

Thus, control and viability theory methods have allowed us to 
introduce ecosystem considerations, such as multispecies and
multiobjectives, and have contributed to integrate the long term dynamics, 
which is generally not considered in conventional fishery management.  
Notice that the World Summit on Sustainable Development
\citep*{Garcia:2003} encouraged the application of the ecosystem approach
by 2010. 

\bigskip

\paragraph*{Acknowledgments.}

This paper was prepared within the MIFIMA (Mathematics, 
  Informatics and Fisheries Management) international research network.
We thank CNRS, INRIA and the French Ministry of Foreign Affairs 
for their funding and support through the regional cooperation program
STIC--AmSud.
We thank the staff of the Peruvian Marine Research Institute (IMARPE),
especially Erich Diaz and Nathaly Vargas for discussions on anchovy and
hake fisheries. 
We also thank Yboon Garcia (IMCA-Peru and CMM-Chile)
for a discussion on the ecosystem model case.

\newcommand{\noopsort}[1]{} \ifx\undefined\allcaps\def\allcaps#1{#1}\fi

 \end{document}